\documentclass[a4paper,11pt]{article}
\usepackage{amsmath,amscd,amssymb,french}

\newcommand{\PP}{{\bf{P}}}
\newcommand{\CC}{{\bf{C}}}
\newcommand{\ZZ}{{\bf{Z}}}

\newcommand{\RR}{{\bf{R}}}

\newcommand{\W}{\textup{W}}
\newcommand{\V}{\textup{V}}
\newcommand{\X}{\textup{X}}
\newcommand{\Y}{\textup{Y}}
\newcommand{\Z}{\textup{Z}}
\newcommand{\K}{\textup{K}}
\newcommand{\R}{\textup{R}}
\newcommand{\C}{\textup{C}}
\newcommand{\N}{\textup{N}}
\newcommand{\F}{\textup{F}}
\newcommand{\U}{\textup{U}}
\newcommand{\E}{\textup{E}}
\newcommand{\G}{\textup{G}}
\newcommand{\I}{\textup{I}}
\newcommand{\A}{\textup{A}}
\newcommand{\D}{\textup{D}}
\newcommand{\T}{\textup{T}}
\newcommand{\M}{\textup{M}}

\newcommand{\re}[2]{(\textit{voir} \cite{#1} #2)}
\newcommand{\res}[1]{(\textit{voir} \cite{#1})}

\newcommand{\qed}{\hfill$\blacksquare$}

\newcounter{compteur}
\newcounter{sect}
\newcommand{\sect}[1]{\\$\ $\newline\stepcounter{sect}\textbf{\thesect}.
\textbf{#1}\setcounter{compteur}{1}\\$\ $\newline}

\newcommand{\theointro}[1]{\\$\ $\newline\textbf{Th\'eor\`eme}.$-$
\noindent\textit{#1}\\$\ $\newline}

\newcommand{\lemme}[1]{\\$\ $\newline\textbf{Lemme
\thesect.\thecompteur}.$-$\noindent\textit{#1}\stepcounter{compteur}\\
$\ $\newline}

\newcommand{\lemmeref}[3]{\\$\ $\newline\textbf{Lemme
\thesect.\thecompteur} (\textit{voir} \cite{#1} #2).$-$
\noindent\textit{#3}\stepcounter{compteur}\\$\ $\newline}

\newcommand{\prop}[1]{\\$\ $\newline\textbf{Proposition
\thesect.\thecompteur}.$-$\noindent\textit{#1}\stepcounter{compteur}\\
$\ $\newline}

\newcommand{\proprefs}[2]{\\$\ $\newline\textbf{Proposition
\thesect.\thecompteur} (\textit{voir} \cite{#1}).$-$
\noindent\textit{#2}\stepcounter{compteur}\\$\ $\newline}

\newcommand{\propenu}[1]{\\$\ $\newline\textbf{Proposition
\thesect.\thecompteur}.$-$\noindent\textit{#1}\stepcounter{compteur}}

\newcommand{\marque}{(\textbf{\thesect.\thecompteur) }\stepcounter{compteur}}

\begin{document}

\centerline{\large{\textbf{Singularit\'es symplectiques}}}
$\ $
\vspace{0.5cm}\\
$\ $
\centerline{St\'ephane \textsc{Druel}}
$\ $
\vspace{1cm}\\
\textbf{Abstract}. We classify isolated symplectic singularities of dimension greater or equal to 6
such that the normalized
blow-up of the singular point is a resolution of singularities whose exceptionnal locus 
is a reduced simple normal crossing divisor with at least two irreducible components. 
There are isomorphic to the quotient singularities of type $\frac{1}{3}(1,2,\ldots,1,2)$.\\
\vspace{1cm}\\
\noindent\textbf{Introduction}\\
\newline
\indent Soit $\textup{(V,\,o)}$ un germe de singularit\'e complexe analytique
\textit{normal}. La
singularit\'e est dite \textit{symplectique} s'il existe une 2-forme
symplectique
$\varphi$ sur le lieu r\'egulier $\V_{\text{reg}}$ de V,
\textit{i.e.} une \textit{2-forme holomorphe ferm\'ee} et \textit{non
  d\'eg\'en\'er\'ee},
et si pour toute r\'esolution $\pi\,:\,\X\longrightarrow\V$ des
singularit\'es de V,
le pull-back de $\varphi$ \`a $\pi^{-1}(\V_{\text{reg}})$ est la
restriction
d'une 2-forme holomorphe sur X \re{Be00}{d\'efinition 1.1}. Si le lieu
singulier de
V est de codimension $\ge 4$ alors la derni\`ere condition est
toujours v\'erifi\'ee \res{Fl88}.\\
\indent Une alg\`ebre de Lie complexe simple a une plus petite orbite
nilpotente non nulle
$\mathcal{O}_{min}$ pour l'action adjointe ; son adh\'erence
$\overline{\mathcal{O}}_{min}=\mathcal{O}_{min}\cup\{0\}$ a une
singularit\'e symplectique
isol\'ee en 0 isomorphe au c\^one sur la vari\'et\'e lisse
$\PP\mathcal{O}_{min}$
des droites de $\mathcal{O}_{min}$ et \textit{toute singularit\'e
symplectique isol\'ee dont le c\^one
tangent est lisse, est analytiquement isomorphe \`a
$(\overline{\mathcal{O}}_{min},\,0)$ pour
une alg\`ebre de Lie complexe simple convenable} \res{Be00}.\\
\indent Fixons un entier $n>0$. Soit $\zeta$ une racine primitive cubique de
l'unit\'e et soit
$<\zeta>$ le groupe cyclique d'ordre 3 engendr\'e par $\zeta$. Soit
$\V=\CC^{2n}\diagup<\zeta>$ o\`u $<\zeta>$ agit sur $\CC^{2n}$ par la
formule\\
\centerline{$\zeta.(z_{1},\,z_{2},\ldots,\,z_{2n-1},\,z_{2n})
=(\zeta z_{1},\,\zeta^{2}z_{2},\ldots,\,\zeta
z_{2n-1},\,\zeta^{2}z_{2n})$.}
La 2-forme symplectique $\sum_{i=1}^{n}dz_{2i-1}\wedge
dz_{2i}$
sur $\CC^{2n}$ est invariante sous $<\zeta>$ et induit une 2-forme
symplectique sur
$\V_{\text{reg}}$. Notons 0 l'unique point singulier de V.
Soit $\pi\,:\,\X\longrightarrow\V$
l'\'eclatement normalis\'e de 0 dans V. Le diviseur
$\pi^{-1}(\text{0})$ est r\'eduit et
globalement \`a croisements normaux (\textit{voir} \S 1). Nous prouvons
le
\theointro{Soit \textup{(V,\,o)} une
singularit\'e symplectique isol\'ee de dimension $2n\ge 6$ et
soit $\pi\,:\,\X\longrightarrow\V$ l'\'eclatement normalis\'e
de \textup{o} dans \V. On suppose que le
diviseur $\pi^{-1}(\textup{o})$ est r\'eduit, globalement \`a
croisements normaux et qu'il a au moins deux composantes
irr\'eductibles. Alors \textup{(V,\,o)} est analytiquement isomorphe 
\`a la singularit\'e quotient $(\CC^{2n}\diagup<\zeta>,\,0)$.}
\indent L'hypoth\`ese $n\ge 3$ permet d'obtenir l'annulation 
de certains groupes de cohomologie (\textit{voir} lemme 3.6).\\
\indent Soit $\F:=\PP_{\PP^{n-1}}(\mathcal{E})$ o\`u
$\mathcal{E}=\mathcal{O}_{\PP^{n-1}}(2)\oplus
\text{H}^{0}(\PP^{n-1},\,\mathcal{O}_{\PP^{n-1}}(1))\otimes
\mathcal{O}_{\PP^{n-1}}$ et soit
$p_{\F}$ le morphisme vers $\PP^{n-1}$.
Soient $p$ et $q$ les projections de $\PP^{n-1}\times\PP^{n-1}$
sur $\PP^{n-1}$ et soit
$i$ l'immersion ferm\'ee $\PP^{n-1}\times\PP^{n-1}\subset\F$
au dessus de $\PP^{n-1}$ correspondant au quotient inversible sur
$\PP^{n-1}\times\PP^{n-1}$\\
\centerline{$\displaystyle{p^{*}\mathcal{E}\twoheadrightarrow
\textup{H}^{0}(\PP^{n-1},\,\mathcal{O}_{\PP^{n-1}}(1))\otimes
\mathcal{O}_{\PP^{n-1}\times\PP^{n-1}}
\twoheadrightarrow q^{*}\mathcal{O}_{\PP^{n-1}}(1)}$}
o\`u $\PP^{n-1}\times\PP^{n-1}$ est
consid\'er\'e comme vari\'et\'e sur
$\PP^{n-1}$ via $p$.
Soit enfin $j:=i\circ s$ o\`u $s$ est l'involution naturelle de
$\PP_{\PP^{n-1}}\times\PP_{\PP^{n-1}}$.
Nous montrons que le lieu exceptionnel est isomorphe au recollement de
deux copies de F
le long de $\PP^{n-1}\times\PP^{n-1}$ via les immersions ferm\'ees $i$ et
$j$ puis nous montrons
que le compl\'et\'e formel de $\X$ le long du diviseur exceptionnel
est d\'etermin\'e \`a isomorphisme pr\`es par le germe de singularit\'e
$(\V,\,\text{o})$. Le r\'esultat cherch\'e est alors une cons\'equence
du premier th\'eor\`eme de comparaison de la th\'eorie ``alg\'ebrique''
\` a la th\'eorie ``formelle'' \re{Gr66}{Chap. III th\'eor\`eme
4.1.5} et du th\'eor\`eme d'approximation d'Artin
\re{Ar68}{corollaire 1.6}.\\
\newline
\noindent\textit{Remerciements}.$-$Je tiens \`a exprimer ici mes
remerciements \`a Arnaud Beauville pour
toute l'aide qu'il m'a apport\'ee.
\sect{Exemple et g\'eom\'etrie torique}
\indent Nous renvoyons \`a \cite{Od88} pour les d\'efinitions et
propri\'et\'es des vari\'et\'es toriques. Fixons un entier $n\ge 1$.
Soit $\N_{0}=\ZZ e_{1}\oplus\cdots\oplus\ZZ e_{2n}$ un $\ZZ$-module
libre de rang $2n$ et soit
$\N=\N_{0}+\frac{\ZZ}{3}(1,\,-1,\ldots,\,1,\,-1)\subset\N_{0}\otimes\RR$.
Soit $\M=\textup{Hom}_{\ZZ}(\N,\ZZ)$. Soit enfin
$\sigma=<e_{1},\ldots,\,e_{2n}>\subset\N\otimes\RR$.
La vari\'et\'e affine $\V$ est torique et son alg\`ebre de fonctions est
$\CC[z_{1},\ldots,\,z_{2n}]^{<\zeta>}=\CC[\M\cap\sigma^{\vee}]$
\re{Od88}{\S 1.5}. Soit $(e_{1}^{*},\ldots,\,e_{2n}^{*})$ la base duale de
$(e_{1},\ldots,\,e_{2n})$. Les fonctions
$$\left\lbrace
\begin{array}{lr}
z_{2i}z_{2j-1} & (i,\,j\in\{1,\ldots,\,n\})\\
z_{2i}z_{2j}z_{2k} & (i,\,j,\,k\in\{1,\ldots,\,n\})\\
z_{2i-1}z_{2j-1}z_{2k-1} & (i,\,j,\,k\in\{1,\ldots,\,n\})
\end{array}
\right.
$$
sont invariantes sous le groupe $<\zeta>$ et forment un syst\`eme
minimal de
g\'en\'erateurs sur $\CC$ de l'alg\`ebre
$\CC[z_{1},\ldots,\,z_{2n}]^{<\zeta>}$, autrement dit, les \'el\'ements
$$\left\lbrace
\begin{array}{lr}
e_{2i}^{*}+e_{2j-1}^{*} & (i,\,j\in\{1,\ldots,\,n\})\\
e_{2i}^{*}+e_{2j}^{*}+e_{2k}^{*} & (i,\,j,\,k\in\{1,\ldots,\,n\})\\
e_{2i-1}^{*}+e_{2j-1}^{*}+e_{2k-1}^{*} & (i,\,j,\,k\in\{1,\ldots,\,n\})
\end{array}
\right.
$$
forment \textit{une base de Hilbert} du monoide $\M\cap\sigma^{\vee}$.
Soient
$e_{2n+1}=\frac{1}{3}(1,\,2,\ldots,\,1,\,2)$ et
$e_{2n+2}=\frac{1}{3}(2,\,1,\ldots,\,2,\,1)$.
L'\'eclatement normalis\'e $\X$ de l'id\'eal $\mathfrak{m}$ de $0$
dans $\V$ est la vari\'et\'e
torique \textit{lisse} dont l'\'eventail est l'ensemble des
c\^ones \textit{r\'eguliers}
$$\left\lbrace
\begin{array}{lr}
\sigma_{i,\,j}=<e_{1},\ldots,\,\widehat{e}_{2i},\ldots,\,\widehat{e}_{2j-1},\ldots,\,e_{2n+2}>
& (i,\,j\in\{1,\ldots,\,n\})\\
\sigma'_{i}=<e_{1},\ldots,\,\widehat{e}_{2i},\ldots,\,\widehat{e}_{2n+1},\,e_{2n+2}>
& (i\in\{1,\ldots,\,n\})\\
\sigma''_{j}=<e_{1},\ldots,\,\widehat{e}_{2j-1},\ldots,\,e_{2n+1},\,\widehat{e}_{2n+2}>
& (j\in\{1,\ldots,\,n\})
\end{array}
\right.
$$
et de leurs faces.
Le lieu exceptionnel est un diviseur globalement \`a croisements
normaux r\'eunion des deux diviseurs irr\'eductibles $\V(e_{2n+1})$ et
$\V(e_{2n+2})$.\\
\indent L'id\'eal $\mathfrak{m}\mathcal{O}_{\X}$ du diviseur
$\pi^{-1}(0)$ est engendr\'e
par la fonction $z_{2i}z_{2j-1}$ (resp. $z_{2i}^3$ et $z_{2j-1}^3$)
sur l'ouvert affine $\U_{\sigma_{2i,\,2j-1}}$ $(i,\,j\in\{1,\ldots,\,n\})$
(resp. $\U_{\sigma'_{i}}$ $(i\in\{1,\ldots,\,n\})$ et
$\U_{\sigma''_{j}}$
$(j\in\{1,\ldots,\,n\})$). Notons $v_{i}=e_{i}$ pour $1\le i\le 2n-1$
et $v_{2n}=e_{2n+2}$.
Soit $(v_{1}^{*},\ldots,v_{2n}^{*})$ la base duale de
$(v_{1},\ldots,v_{2n})$. Les \'el\'ements
$(v_{i})_{1\le i\le 2n}$ forment une base du
$\ZZ$-module $\M$ et l'alg\`ebre des fonctions r\'eguli\`eres sur
l'ouvert $\U_{\sigma'_{n}}$ est l'alg\`ebre sur $\CC$ du monoide libre
engendr\'e par ces \'el\'ements.
On a $v_{2i}^{*}=e_{2i}^{*}-e_{2n}^{*}$ pour $1\le i\le n-1$,
$v_{2i-1}^{*}=e_{2i-1}^{*}-2e_{2n}^{*}$ pour $2\le i\le n$ et
$v_{2n}^{*}=3e_{2n}^{*}$. Le diviseur $\pi^{-1}(0)$ est donc
localement d\'efini par l'annulation d'une coordonn\'ee sur l'ouvert
consid\'er\'e et en particulier r\'eduit le long de $\V(e_{2n+2})$. On
v\'erifie que ce diviseur est
\'egalement r\'eduit le long de $\V(e_{2n+1})$.
\sect{Pr\'eliminaires}
\marque Soit X une vari\'et\'e projective lisse sur le corps
$\CC$ des nombres complexes. Soit
$\N_{1}(\X)=(\{\text{1-cycles}\}/\equiv)\otimes\RR$
\noindent o\`u $\equiv$ d\'esigne l'\'equivalence num\'erique.
On consid\`ere le c\^one $\textup{NE}(\X)\subset \N_{1}(\X)$
engendr\'e par les classes des 1-cycles effectifs. Une \textit{ar\^ete
extr\'emale} est une demi-droite $\R$ dans $\overline{\textup{NE}}(\X)$,
adh\'erence de $\textup{NE}(\X)$ dans $\N_{1}(\X)$, v\'erifiant
$\K_{\X}.\R^{*}<0$  et telle que pour tout
$\Z_{1},\Z_{2}\in\overline{\textup{NE}}(\X)$, si
$\Z_{1}+\Z_{2}\in \R$ alors $\Z_{1},\,\Z_{2}\in \R$.
Une \textit{courbe rationnelle extr\'emale} est une
courbe rationnelle irr\'eductible $\C$ telle que
$\RR^{+}[\C]$ soit une ar\^ete extr\'emale et
$-\K_{\X}.\C\le\text{dim}(\X)+1$. \textit{Toute
ar\^ete extr\'emale} R
\textit{est engendr\'ee par une courbe rationnelle
extr\'emale} et admet une
\textit{contraction}, c'est-\`a-dire qu'il existe une
vari\'et\'e projective normale $\Y$ et un morphisme
$\phi\,:\,\X\longrightarrow\Y$, surjectif \`a fibres connexes,
contractant les courbes irr\'eductibles $\C$ telles que $[\C]\in\R$
(th\'eor\`eme de Kawamata-Shokurov).\\
\indent Rappelons un r\'esultat de
J. Wisniewski \res{Wi91} sur le
lieu exceptionnel d'une contraction extr\'emale. Soit $\F$ une
composante irr\'eductible d'une fibre non triviale d'une contraction
\'el\'ementaire associ\'ee \`a l'ar\^ete extr\'emale $\R$. Nous appelons
\textit{lieu de $\R$}, le lieu des courbes dont la classe
d'\'equivalence
num\'erique appartient \`a $\R$. On a alors l'in\'egalit\'e :
$$\text{dim}(\F)+\text{dim}(\text{lieu de
}\R)\ge\text{dim}(\X)+\ell(\R)-1,$$
\noindent o\`u $\ell(\R)$ d\'esigne la \textit{longueur} de l'ar\^ete
extr\'emale R :
$$\ell(\R)=\text{inf}\{-\K_{\X}.\C_{0}\,|\,\C_{0} \text{ \'etant une
courbe
rationnelle et }\C_{0}\in \R\}.$$
\lemme{Soit $\X$ une vari\'et\'e projective de dimension
$n\ge 2$ et soit $\Y\subset\X$ un diviseur effectif. Soit $\textup{D}$
un diviseur
num\'eriquement effectif. Si
$-\K_{\X}\equiv\textup{D}+\Y$ alors il existe une ar\^ete extr\'emale
$\R$ telle que $\Y.\R^{*}>0$.}
\textit{D\'emonstration}.$-$Soit $\C\subset\X$ une courbe telle que
$\Y.\C>0$. La courbe $\C$ se d\'ecompose
$\C\equiv\C_{0}+\sum_{i\in\I}a_{i}\C_{i}$ avec
$\C_{0}\in\overline{\text{NE}}^{+}(\X)=
\{z\in\overline{\text{NE}}(\X)\,|\,\K_{\X}.z\ge0\}$ et
$(\C_{i})_{i\in\I}$ est l'ensemble des courbes rationelles extr\'emales
\re{Mo82}{th\'eor\`eme 1.5}. Les coefficients $(a_{i})_{i\in\I}$ sont
positifs ou nuls et presque tous nuls. On a
$\Y.\C_{0}=-\K_{\X}.\C_{0}-\textup{D}.\C_{0}\le 0$ et il existe donc un
\'el\'ement $i\in\I$ tel que $a_{i}\C_{i}.\Y>0$.\qed\\
\newline
\marque Fixons quelques notations. Soient $a\in\ZZ$ et $k\ge 2$. Soit
$\F_{a}:=\PP_{\PP^{k-1}}(\mathcal{E}_{a})$ o\`u
$\mathcal{E}_{a}=\mathcal{O}_{\PP^{k-1}}(a)\oplus
\text{H}^{0}(\PP^{k-1},\,\mathcal{O}_{\PP^{k-1}}(1))\otimes\mathcal{O}_{\PP^{k-1}}$ et soit
$p_{a}$ le morphisme vers $\PP^{k-1}$.
Soient $p$ et $q$ les projections de $\PP^{k-1}\times\PP^{k-1}$
sur $\PP^{k-1}$ et soit enfin $i_{a}$ l'immersion ferm\'ee
$\PP^{k-1}\times\PP^{k-1}\subset\F_{a}$
au dessus de $\PP^{k-1}$ correspondant au quotient inversible sur
$\PP^{k-1}\times\PP^{k-1}$\\
\centerline{$\displaystyle{p^{*}\mathcal{E}_{a}\twoheadrightarrow
\textup{H}^{0}(\PP^{k-1},\,\mathcal{O}_{\PP^{k-1}}(1))\otimes
\mathcal{O}_{\PP^{k-1}\times\PP^{k-1}}
\twoheadrightarrow q^{*}\mathcal{O}_{\PP^{k-1}}(1)}$}
o\`u $\PP^{k-1}\times\PP^{k-1}$ est
consid\'er\'e comme vari\'et\'e sur
$\PP^{k-1}$ via $p$.
\propenu{Soit $(\W,\,\textup{L})$ une vari\'et\'e polaris\'ee de dimension
$m\ge 2$ et soit
$\Y=\cup_{i=1}^{l}\Y_{i}$ un diviseur r\'eduit et globalement \`a croisements
normaux.
On suppose $\omega_{\W}\simeq\textup{L}^{-k}(-\Y)$.
\begin{enumerate}
\item Si $k>\frac{m+1}{2}$ et $l=1$ ou
$k\ge\frac{m+1}{2}$ et $l\ge 2$ alors $\W$ est de Fano et
$b_{2}(\W)=1$.
\item Si $k=\frac{m+1}{2}$ et $l=1$ alors $\Y\subset \W$ est
isomophe \`a $\textup{Q}_{2}\subset\PP^{3}$ ou
$\Y\subset \W$ est isomophe \`a $\textup{Q}_{2}\subset\textup{Q}_{3}$ ou
$\Y$ et $\W$ sont de Fano et $b_{2}(\Y)=b_{2}(\W)=1$ ou
$\Y\subset\W$ est isomorphe
\`a $\PP^{k-1}\times\PP^{k-1}\subset\F_{a}$ pour un entier $a\ge 0$
convenable,  $\textup{L}$ est
isomorphe au fibr\'e $\mathcal{O}_{\F_{a}}(1)\otimes
p_{a}^{*}(\mathcal{O}_{\PP^{k-1}}(1))$
et l'id\'eal de $\Y$ dans $\W$ est $\textup{L}^{\otimes -1}\otimes
p_{a}^{*}\mathcal{O}_{\PP^{k-1}}(a+1)$.
\end{enumerate}
}
\noindent\textit{D\'emonstration}.$-$La premi\`ere partie de la
proposition se montre par r\'ecurrence sur la
dimension $m\ge 2$ de $\W$. Si $m=2$ alors
$k\ge 2$ et, par le lemme 2.2, il existe une ar\^ete extr\'emale $\R$
telle que $\R^{*}.\Y>0$ et
donc de longueur $\ge 3$. La surface $\W$ est donc de Fano et de nombre
de Picard $1$
\re{Wi89}{proposition 2.4.1}.\\
\indent Supposons la proposition d\'emontr\'ee jusqu'en dimension
$m-1\ge 2$ et consid\'erons
$\W$ comme dans l'\'enonc\'e. Soit $\C_{1}$ une courbe rationnelle
extr\'emale telle que $\Y.\C_{1}>0$ et
$\R_{1}=\RR_{+}[\C_{1}]$ et supposons par exemple $\Y_{1}.\C_{1}>0$
(\textit{voir} lemme 2.2).
Notons $\W\overset{\varphi_{1}}{\longrightarrow}\Z_{1}$ la contraction
\'el\'ementaire correspondante.
La formule d'adjonction donne
$\omega_{\Y_{1}}\simeq\textup{L}_{|\Y_{1}}^{-k}(-\sum_{i\neq 1}\Y_{i}\cap\Y_{1})$
avec
$k\ge\frac{m+1}{2}>\frac{\text{dim}(\Y_{1})+1}{2}$. Si $\sum_{i\neq
1}\Y_{i}\cap\Y_{1}\neq\emptyset$ alors,
par hypoth\`ese de r\'ecurrence, $\Y_{1}$ est de Fano et
$b_{2}(\Y_{1})=1$. Si
$\sum_{i\neq 1}\Y_{i}\cap\Y_{1}=\emptyset$ alors $\Y_{1}$ est de Fano
et $b_{2}(\Y_{1})=1$ sauf si
$k=\frac{m+1}{2}$ et $(\Y_{1},\,\textup{L}_{|\Y_{1}})
\simeq(\PP^{k-1}\times\PP^{k-1},
\,\mathcal{O}_{\PP^{k-1}}(1)\boxtimes\mathcal{O}_{\PP^{k-1}}(1))$
\res{Wi90}.\\
\indent Soit $\F_{1}$ la fibre de $\varphi_{1}$ contenant $\C_{1}$.
La fibre $\F_{1}$ rencontre en particulier $\Y_{1}$ et
$\text{dim}(\F_{1}\cap\Y_{1})\ge\text{dim}(\F_{1})-1\ge\ell(\R_{1})-2\ge
\frac{m-1}{2}>0$.
Soit $\C\subset\F_{1}\cap\Y_{1}$ une courbe irr\'eductible. La courbe
$\C$ est contract\'ee
par $\varphi_{1}$ et on a donc $\C\in\R_{1}$ et $\C.\Y_{1}>0$. Notons
que
par tout point de $\Y_{1}$ il passe une courbe trac\'ee sur $\Y_{1}$
num\'eriquement
proportionnelle \`a $\C$ et donc contract\'ee par $\varphi_{1}$. La
contraction
\'el\'ementaire $\varphi_{1}$ est donc divisorielle ou de type fibr\'ee.
Si $\varphi_{1}$ est divisorielle
alors $\R_{1}$ n'est pas num\'eriquement effective, autrement dit, il
existe un diviseur
irr\'eductible $\D$ tel que $\D.\C<0$. On a donc $\C\subset\D$ et
$\Y_{1}=\D$, ce qui est impossible
puisque $\C.\Y_{1}>0$. La contraction $\varphi_{1}$ est donc de type
fibr\'ee.\\
\noindent\textit{Supposons $b_{2}(\Y_{1})=1$}.$-$Le diviseur $\Y_{1}$
est donc contract\'e sur un point
par $\varphi_{1}$ et, puisque $\Y_{1}.\C>0$, $\text{dim}(\Z_{1})=0$. La
vari\'et\'e X est donc de Fano de nombre de
Picard 1.\\
\noindent\textit{Supposons $\Y_{1}\simeq\PP^{k-1}\times\PP^{k-1}$ et
$\textup{L}_{|\Y_{1}}\simeq\mathcal{O}_{\PP^{k-1}}(1)\boxtimes
\mathcal{O}_{\PP^{k-1}}(1)$}.$-$Notons
$(a,b)$ le bidegr\'e de la courbe $\C\subset\Y_{1}$.\\
\indent Supposons $a>0$ et $b>0$. Alors, par tout couple de points de $\Y_{1}$, il
passe une courbe
trac\'ee sur $\Y_{1}$ de bidegr\'e $(a,b)$ et donc num\'eriquement
proportionnelle \`a $\C$ dans $\W$. Le diviseur $\Y_{1}$ est donc
contract\'e sur un
point par $\varphi_{1}$. Le morphisme
$\Y_{1}\overset{{\varphi_{1}}_{|\Y_{1}}}{\longrightarrow}\Z_{1}$ est
surjectif puisque
$\Y_{1}.\C>0$. Finalement, $\Z_{1}$ est de dimension 0 et
$b_{2}(\W)=1$. On en
d\'eduit que $k=2$ par le th\'eor\`eme de Lefschetz, puis que
$\W$ est une quadrique  et $\Y_{1}$ une section hyperplane ou que $\W$
est un espace
projectif et $Y_{1}$ une quadrique par le th\'eor\`eme de
Kobayashi-Ochiai \res{KO73}.\\
\indent Supposons maintenant $(a,b)=(0,1)$. Notons $p$ la premi\`ere
projection de $\Y_{1}$ sur $\PP^{k-1}$. Soit $\F\subset\Y_{1}$ une fibre
ensembliste
de ${\varphi_{1}}_{|\Y_{1}}$. Les courbes trac\'ees sur $\Y_{1}$ de
bidegr\'e
$(0,1)$ sont contract\'ees par $\varphi_{1}$ et on a
donc $\F=p^{-1}(p(\F))$. S'il existe une courbe $\C_{2}$ trac\'ee sur
$\F$ de bidegr\'e $(a_{1},b_{1})$ avec
$a_{1}>0$ alors le 1-cycle $\C_{2}+\C_{1}\in\R_{1}$ est de bidegr\'e
$(a_{1},b_{1}+1)$ avec $a_{1}>0$ et $b_{1}+1>0$
et les arguments utilis\'es ci-dessus suffisent pour conclure. On peut
donc supposer que les
courbes trac\'ees sur $\F$ sont de bidegr\'e $(0,*)$. On en d\'eduit
que $p(\F)$ est de dimension 0 puis que $\F$ est \'equidimensionnelle
de dimension $k-1$ et que $\Z_{1}$ est de dimension $k-1$.
Soit $\F_{1}\subset\W$ une composante irr\'eductible d'une
fibre de $\varphi_{1}$ ; $\F_{1}$ rencontre $\Y_{1}$ puisque
$\R^{*}.\Y_{1}>0$ et $\text{dim}(\F_{1}\cap\Y_{1})=k-1$. On a donc
$\text{dim}(\F_{1})=k$. Le morphisme $\varphi_{1}$ est
donc \'equidimentionnel de dimension $k$. Soit
$\F_{1}$ une fibre g\'en\'erique de $\varphi_{1}$.
On a
$\omega_{\F_{1}}=\textup{L}_{|\F_{1}}^{-k}\otimes\mathcal{O}_{\F_{1}}(-\Y))$.
Le fibr\'e $\mathcal{O}_{\F_{1}}(\Y))$ est num\'eriquement effectif
puisque
$\R_{1}^{*}.\Y>0$ et $\F_{1}$ est donc de Fano. Soit
$\R_{\F_{1}}$ une ar\^ete extr\'emale de $\F_{1}$. On a
$\ell(\R_{\F_{1}})\ge\text{dim}(\F_{1})+1$ et donc $b_{2}(\F_{1})=1$
\re{Wi89}{proposition 2.4.1}. On en d\'eduit
$(\F_{1},\,\textup{L}_{|\F_{1}})=(\PP^{k},\,\mathcal{O}_{\PP^{k}}(1))$
puis
que $\varphi_{1}$ est un fibr\'e en espaces projectifs
localement trivial pour la topologie de Zariski \re{Fu87}{lemme 2.12}
et en particulier que $\Z_{1}$ est lisse. La fibre ensembliste
$\F_{1}\subset\PP^{k}$ de
${\varphi_{1}}_{|\Y_{1}}$ est r\'eunion disjointes de fibres de $p$ de
dimension $k-1$ et
donc irr\'eductible. Le morphisme
${\varphi_{1}}_{|\Y_{1}}$ s'identifie finalement \`a la seconde
projection
de $\Y_{1}=\PP^{k-1}\times\PP^{k-1}$ sur $\PP^{k-1}$ et $\W$ \`a la
fibration $\PP_{\PP^{k-1}}({\varphi_{1}}_{*}\textup{L})$
au dessus de $\PP^{k-1}$. L'immersion ferm\'ee $\Y_{1}\subset\W$ au
dessus
de $\PP^{k-1}$ est donn\'ee par le quotient
inversible $p^*({\varphi_{1}}_{*}\textup{L})\twoheadrightarrow
\textup{L}_{|\Y_{1}}=\mathcal{O}_{\PP^{k-1}}(1)\boxtimes
\mathcal{O}_{\PP^{k-1}}(1)$ et on a donc
un morphisme surjectif ${\varphi_{1}}_{*}\textup{L}\twoheadrightarrow
p_{*}(\textup{L}_{|\Y_{1}})=\textup{H}^{0}(\PP^{k-1},\,\mathcal{O}_{\PP^{k-1}}(1))
\otimes\mathcal{O}_{\PP^{k-1}}(1)$.
Finalement
${\varphi_{1}}_{*}\textup{L}=\mathcal{O}_{\PP^{k-1}}(a)\oplus
\textup{H}^{0}(\PP^{k-1},\,\mathcal{O}_{\PP^{k-1}}(1))\otimes\mathcal{O}_{\PP^{k-1}}(1)$ avec $a\ge 1$.
L'id\'eal de $\Y_{1}$ dans $\W$ est
$\textup{L}^{\otimes -1}\otimes p^{*}\mathcal{O}_{\PP^{k-1}}(b)$ o\`u
$b$ est un entier convenable.
Le fibr\'e canonique est donn\'e par la formule
$\omega_{\W}=\textup{L}^{-k-1}\otimes
p^{*}\mathcal{O}_{\PP^{k-1}}(a)=
\textup{L}^{-k}(-\Y)$ et l'id\'eal de $\Y$ dans $\W$ est donc
$\textup{L}^{\otimes -1}\otimes
p^{*}\mathcal{O}_{\PP^{k-1}}(a)$. L'id\'eal du diviseur $\sum_{i\neq
1}\Y_{i}$ dans $\W$ est donc
$p^{*}\mathcal{O}_{\PP^{k-1}}(a-b)$. Or les diviseurs $\Y_{1}$ et
$\sum_{i\neq 1}\Y_{i}$ sont
disjoints et on a donc $a=b$ et $\sum_{i\neq 1}\Y_{i}=\emptyset$. La
premi\`ere partie de
la proposition est donc d\'emontr\'ee.\\
\indent La seconde partie de la proposition a en fait d\'ej\`a \'et\'e
prouv\'ee ci-dessus.\qed
\sect{D\'emonstration du th\'eor\`eme}
\marque Fixons les notations. Soit $\textup{(V,\,o)}$ une singularit\'e
symplectique isol\'ee de
dimension $2n$ ($n>1$) v\'erifiant les hypoth\`eses du th\'eor\`eme et
soit
$\pi\,:\,\X\longrightarrow\V$ l'\'eclatement normalis\'e de l'id\'eal
maximal de o dans $\V$. On peut toujours supposer qu'il existe une
2-forme
symplectique $\varphi$ sur $\V-\{o\}$ restriction d'une 2-forme
holomorphe sur $\X$ not\'ee encore $\varphi$.
Notons $\E=\E_{1}\cup\cdots\cup\E_{l}$ $(l\ge 2)$ le diviseur
exceptionnel et $\textup{L}:=\mathcal{O}_{\X}(-\E)$.
Le fibr\'e $\textup{L}_{|\E}$ est ample et engendr\'e par ses sections.
\proprefs{Be00}{L'ordre d'annulation de $\varphi^{n}$ le long de
$\E_{i}$
est $n-1$ pour $1\le i\le l$.}
\textit{D\'emonstration}.$-$Notons $k_{i}$ $(1\le i\le l)$ l'ordre
d'annulation de $\varphi^{n}$
le long de $\E_{i}$ $(1\le i\le l)$ de sorte que
$\text{div}(\varphi^n)=\sum_{i=1}^{l}k_{i}\E_{i}$.\\
\indent Les singularit\'es symplectiques sont rationnelles
\re{Be00}{proposition 1.3}
et ainsi \\ 
$\varphi\in\textup{H}^{0}(\X,\Omega_{\X}^{2}(\text{log E})(-\E))$
\re{CF01}{th\'eor\`eme 2.1} o\`u
$\Omega_{\X}^{2}(\text{log E})$ est le faisceau des 2-formes
diff\'erentielles
m\'eromorphe \`a p\^oles au pire logarithmiques le long de
E. La 2n-forme\\
$\varphi^{n}\in\textup{H}^{0}(\X,\,\omega_{\X}(-(n-1)\E)$ est
non nulle en dehors de $\E$ et il s'ensuit
$k_{i}\ge n-1$ $(1\le i\le l)$ puisque le diviseur E est contract\'e
par $\pi$.\\
\indent Le produit exterieur avec $\varphi^{n-1}$ donne une application
$\mathcal{O}_{\X}$-lin\'eaire
$\Omega_{\X}^{1}\longrightarrow\T_{\X}(\sum_{i=1}^{l}k_{i}\E_{i})$
qui est un isomorphisme en dehors de $\E$.
Notons $k_{i}-j_{i}\ge 0$ l'ordre
d'annulation de cette application le long du diviseur $\E_{i}$ $(1\le
i\le l)$. On obtient une application
$\lambda\,:\,\Omega_{\X}^{1}\longrightarrow\T_{\X}(\sum_{i=1}j_{i}\E_{i})$
dont la restriction au diviseur $\E_{i}$ est non nulle et une section
$\text{det}(\lambda)$ de
$\mathcal{O}_{\X}(2\sum_{i=1}^{l}(nj_{i}-k_{i})\E_{i})$ non nulle
en dehors de $\E$. On a donc $k_{i}\le nj_{i}$ pour $1\le i\le l$
puisque le diviseur $\E$ est contract\'e par $\pi$ et en particulier
$j_{i}\ge 1$.\\
\indent Supposons par exemple $j_{1}\ge j_{i}$ pour $1\le i\le l$ et
consid\'erons le diagramme de
suites exactes\\
$\ $
\newline
\begin{equation*}
\begin{CD}
0 @))) \mathcal{O}_{\E_{1}}(-\E_{1}) @))) {\Omega_{\X}^{1}}_{|\E_{1}}
@))) \Omega_{\E_{1}}^{1} @))) 0 \\
 & & & & @VV{{\lambda}_{|\E_{1}}}V \\
0 @))) \T_{\E_{1}}(\displaystyle{\sum_{i=1}^{l}j_{i}\E_{i}}) @)))
{\T_{\X}}_{|\E_{1}}(\displaystyle{\sum_{i=1}^{l}j_{i}\E_{i}}) @)))
\mathcal{O}_{\E_{1}}((j_{1}+1)\E_{1}+\displaystyle{\sum_{i=2}^{l}j_{i}\E_{i}}) @)))
0. \\
\end{CD}
\end{equation*}\\
$\ $
\newline
On a
$\textup{Hom}_{\E_{1}}(\mathcal{O}_{\E_{1}}(-\E_{1}),
\,\mathcal{O}_{\E_{1}}((j_{1}+1)\E_{1}+\sum_{i=2}^{l}j_{i}\E_{i}))
\subset\textup{H}^{0}(\E_{1},\,\textup{L}_{|\E_{1}}^{-j_{1}-2})=0$ et, de m\^eme,
$\textup{Hom}_{\E_{1}}(\mathcal{O}_{\E_{1}}(-\E_{1}),\T_{\E_{1}}(\sum_{i=1}^{l}j_{i}\E_{i}))=
\textup{Hom}_{\E_{i}}(\Omega_{\E_{1}}^{1},
\,\mathcal{O}_{\E_{1}}((j_{1}+1)\E_{1}+\sum_{i=2}^{l}j_{i}\E_{i}))\subset
\textup{H}^{0}(\E_{1},\T_{\E_{1}}\otimes\textup{L}_{|\E_{1}}^{-j_{1}-1})=0$ puisque
$j_{1}\ge 1$
\re{Be00}{lemme 3.3}.
L'application $\lambda_{|\E_{1}}$ se factorise \`a travers une fl\`eche
antisym\'etrique
$\Omega_{\E_{1}}^{1}\overset{\mu_{\E_{1}}}{\longrightarrow}\T_{\E_{1}}
(\sum_{i=1}^{l}j_{i}\E_{i})$
et provient donc d'un \'el\'ement non nul de
$\textup{H}^{0}(\E_{1},\overset{2}{\wedge}\T_{\E_{1}}(\sum_{i=1}^{l}j_{i}\E_{i}))
=\textup{H}^{0}(\E_{1},\overset{2}{\wedge}\T_{\E_{1}}\otimes\textup{L}_{|\E_{1}}^{-j_{1}}
(-\sum_{i=2}^{l}(j_{1}-j_{i})\E_{i}))$.
Il s'ensuit $j_{1}\le 2$ \re{Be00}{lemme 3.3}. Le diviseur E est connexe
et
l'ensemble $\{\E_{i_{1}},\ldots,\E_{i_{m}}\}$ des diviseurs rencontrant
$\E_{1}$ et distincts de
celui-ci est donc non vide. Si $j_{1}=2$ alors
$j_{i_{1}}=\cdots=j_{i_{m}}=j_{1}=2$,
$(\E_{1},\,\textup{L}_{|\E_{1}})\simeq
(\PP^{2n-1},\,\mathcal{O}_{\PP^{2n-1}}(1))$ et
$(\E_{i_{r}},\,\textup{L}_{|\E_{i_{r}}})\simeq
(\PP^{2n-1},\,\mathcal{O}_{\PP^{2n-1}}(1))$
pour $1\le r\le m$ \re{Be00}{lemme 3.3}. Les fibr\'es
$\mathcal{O}_{\E_{1}}(-\E_{1})\simeq\textup{L}_{|\E_{1}}
(\sum_{i=2}^{l}\E_{i})$ et
$\mathcal{O}_{\E_{i_{1}}}(\E_{1})$ sont en particulier amples ce qui est
absurde
puisque $\E_{1}\cap\E_{i_{1}}$ est de dimension $\ge 1$.
Finalement $j_{i}=1$ et $k_{i}\le n$ pour $1\le i\le l$.
Si $k_{1}=n$ alors $\text{det}(\lambda_{|\E_{1}})$ est
g\'en\'eriquement non nul puisque le
diviseur $\E$ est contract\'e par
$\pi$. Ceci est \`a nouveau exclu car $\lambda_{|\E_{1}}$
s'annule sur le sous-fibr\'e
$\mathcal{O}_{\E_{1}}(-\E_{1})\subset{\Omega_{\X}^{1}}_{|\E_{1}}$.
Ceci d\'emontre les \'egalit\'es $k_{i}=n-1$ pour $1\le i\le l$.\qed
\lemme{Le lieu exceptionnel $\E$ de $\pi$ a exactement deux composantes
irr\'eductibles dont l'intersection
est connexe.}
\textit{D\'emonstration}.$-$Le fibr\'e canonique $\omega_{\X}$ de $\X$
est isomorphe \`a
$\textup{L}^{-(n-1)}$ par la proposition 3.2 et
$\omega_{\E_{i}}\simeq\textup{L}_{|\E_{i}}^{-n}(-\sum_{j\neq
i}\E_{i}\cap\E_{i})$
par la formule d'adjonction. Si $l\ge 3$ alors il existe
$i_{1}\in\{1,\ldots,l\}$ tel que
$\E_{i_{1}}$ rencontre au moins deux autres composantes irr\'eductibles
de $\E$ et
$\E_{i_{1}}$ est de Fano avec $b_{2}(\E_{i_{1}})=1$ par la proposition
2.4. Notons $\{\E_{i_{1}},\ldots,\E_{i_{m}}\}$ l'ensemble
des composantes rencontrant $\E_{i_{1}}$ $(m\ge 3)$. Le diviseur 
$\E_{i_{1}}\cap\E_{i_{2}}\subset\E_{i_{1}}$ est ample et en particulier
connexe. Si $b_{2}(\E_{i_{2}})=1$ alors
$\E_{i_{2}}$ est de Fano et le fibr\'e
$\mathcal{O}_{X}(\E_{i_{1}})_{|\E_{i_{2}}}\simeq
\mathcal{O}_{\E_{i_{2}}}(\E_{i_{1}})$
et sa restriction \`a $\E_{i_{1}}\cap\E_{i_{2}}$ sont donc amples.
Or le fibr\'e 
$\mathcal{O}_{X}(-\E_{i_{1}})_{|\E_{i_{1}}}\simeq\mathcal{O}_{\E_{i_{1}}}(-\E_{i_{1}})$
est isomorphe au fibr\'e $\textup{L}_{\E_{i_{1}}}(\sum_{j\neq i_{1}}\E_{j})$ et sa restriction
\`a $\E_{i_{1}}\cap\E_{i_{2}}$ est \'egalement ample, ce qui est absurde puisque 
$\E_{i_{1}}\cap\E_{i_{2}}$ est de dimension $\ge 1$.
Il existe donc un isomorphisme
$(\E_{i_{1}}\cap\E_{i_{2}}\subset\E_{i_{2}})\simeq
(\PP^{n-1}\times\PP^{n-1}\subset\PP_{\PP^{n-1}}(\mathcal{E}_{a})$ pour un
entier $a\ge 0$ convenable
(\textit{voir} proposition 2.4). Le fibr\'e
$\mathcal{O}_{\E_{i_{2}}}(\E_{i_{1}})$
est isomorphe au fibr\'e $\textup{L}\otimes
p_{a}^{*}\mathcal{O}_{\PP^{k-1}}(-a-1)$
(\textit{voir} proposition 2.4) et sa restriction \`a
$\E_{i_{1}}\cap\E_{i_{2}}\simeq\PP^{n-1}\times\PP^{n-1}$ est isomorphe 
au fibr\'e $\mathcal{O}_{\PP^{n-1}}(1)\boxtimes\mathcal{O}_{\PP^{n-1}}(-a)$
qui n'est pas ample, ce qui est \`a nouveau absurde.
La connexit\'e de l'intersection des deux composantes irr\'eductibles de
$\E$ se d\'emontre par des
arguments analogues.\qed\\
\newline
\marque Posons $i=i_{2}$ et $\mathcal{E}=\mathcal{E}_{2}$. Soit enfin
$j:=i\circ s$ o\`u $s$ est l'involution naturelle de
$\PP^{n-1}\times\PP^{n-1}$. Notons $\F$ et $\G$ les deux composantes
irr\'eductibles de $\E$. Posons $\textup{W}:=
\textup{H}^{0}(\PP^{n-1},\,\mathcal{O}_{\PP^{n-1}}(1))$.
\prop{Le lieu exceptionnel $\E$ de $\pi$ est isomorphe au recollement
de deux copies de $\PP_{\PP^{n-1}}(\mathcal{E})$ le long de
$\PP^{n-1}\times\PP^{n-1}$ via les immersions ferm\'ees $i$ et $j$.}
\textit{D\'emonstration}.$-$Notons $\textup{H}=\F\cap\G$. La
proposition 3.2 et la formule
d'adjonction donnent
$\omega_{\F}\simeq\textup{L}_{|\F}^{-n}(-\textup{H})$ et
$\omega_{\G}\simeq\textup{L}_{|\G}^{-n}(-\textup{H})$. Si
$b_{2}(\textup{H})=1$ alors
$\F$ et $\G$ sont de Fano et $b_{2}(\F)=b_{2}(\G)=1$ par la
proposition 2.4.
Les fibr\'es $\textup{L}_{|\F}$ et $\mathcal{O}_{\F}(\G)$ sont amples et
$\mathcal{O}_{\F}(-\F)$ l'est
donc \'egalement. Le fibr\'e $\mathcal{O}_{\G}(\F)$ est ample et
$\textup{H}=\F\cap\G$ est de dimension
au moins 1, ce qui est finalement absurde. On a donc
$b_{2}(\textup{H})\ge 2$ et H est
isomorphe \`a $\PP^{n-1}\times\PP^{n-1}$ toujours par la proposition
2.4. Supposons $\F$ isomorphe \`a $\PP^{3}$ ou $\textup{Q}_{3}$
($n=2$). La restriction du fibr\'e
$\mathcal{O}_{\X}(-\F)=\textup{L}\otimes\mathcal{O}_{\X}(\G)$ \`a
$H=\PP^{1}\times\PP^{1}$ est
en particulier ample et on a donc $b_{2}(\G)\ge 2$, puisque
$\mathcal{O}_{\X}(\F)_{|\G}$ est effectif et $\F\cap\G$ est de dimension
$\ge 1$,
et $\textup{H}\subset\G$ est
isomorphe \`a
$\PP^{1}\times\PP^{1}\subset\PP_{\PP^{1}}(\mathcal{E}_{a_{\G}})$.
La restriction de $\mathcal{O}_{\X}(-\F)$
\`a $\PP^{1}\times\PP^{1}$ est
$\mathcal{O}_{\PP^{1}}(a_{\G})\boxtimes\mathcal{O}_{\PP^{1}}(-1)$
ou $\mathcal{O}_{\PP^{1}}(-1)\boxtimes\mathcal{O}_{\PP^{1}}(a_{\G}))$
(\textit{voir} proposition 2.4)
qui n'est pas ample, ce qui est absurde.\\
\indent Il existe ainsi un isomorphisme de $\textup{H}\subset\F$ (resp.
$\textup{H}\subset\G$) sur
$\PP^{n-1}\times\PP^{n-1}\subset\PP_{\PP^{n-1}}(\mathcal{E}_{a_{\F}})$
(resp.
$\PP^{n-1}\times\PP^{n-1}\subset\PP_{\PP^{n-1}}(\mathcal{E}_{a_{\G}})$
tel que la restriction de $\textup{L}$ \`a F (resp. G) soit isomorphe
\`a
$\mathcal{O}_{\F}(1)\otimes p_{\F}^{*}\mathcal{O}_{\PP^{n-1}}(1)$
(resp. $\mathcal{O}_{\G}(1)\otimes
p_{\G}^{*}\mathcal{O}_{\PP^{n-1}}(1)$) o\`u
$p_{\F}$ (resp. $p_{\G}$) est la projection de $\F$ (resp. $\G$) sur
$\PP^{n-1}$. Enfin,
$\mathcal{O}_{\F}(-\text{H})\simeq\textup{L}_{|\F}^{-1}\otimes
p_{\F}^{*}(\mathcal{O}_{\PP^{n-1}}(a_{\F}+1))$ et
$\mathcal{O}_{\G}(-\textup{H})\simeq\textup{L}_{|\G}^{-1}\otimes
p_{\G}^{*}(\mathcal{O}_{\PP^{n-1}}(a_{\G}+1))$.\\
\indent Soit $\ell\subset\textup{H}$ une droite, \textit{i.e.}
$(\ell.\textup{L}):=\text{deg}_{\ell}({\textup{L}}_{|\ell})=1$,
verticale pour $p_{\G}$. Le fibr\'e
${\mathcal{O}_{\X}(\F)}_{|\G}\simeq\mathcal{O}_{\G}(\textup{H})$ est
isomorphe
\`a $\textup{L}_{|\G}\otimes
p_{\G}^{*}(\mathcal{O}_{\PP^{n-1}}(-a_{\G}-1))$
et en particulier $(\ell.\mathcal{O}_{\X}(\F))=1$. La restriction de
$\mathcal{O}_{\X}(\F)$ \`a $\F$ est isomorphe \`a
$\textup{L}_{|\F}^{-2}\otimes
p_{\F}^{*}(\mathcal{O}_{\PP^{n-1}}(a_{\F}+1))$
et donc
$(\ell.\mathcal{O}_{\X}(\F))=-2+(a_{\F}+1)\text{deg}_{\ell}
(p_{\F}^{*}(\mathcal{O}_{\PP^{n-1}}(1)_{|\ell})$. La droite
$\ell$ est donc horizontale pour $p_{\F}$ et la restriction de
$p_{\F}\times p_{\G}$ \`a H est un isomorphisme de H sur
$\PP^{n-1}\times\PP^{n-1}$, autrement dit, les restrictions de $p_{\F}$
et $p_{\G}$ \`a H s'identifie aux deux projections de
$\PP^{n-1}\times\PP^{n-1}$ sur $\PP^{n-1}$ via les deux isomorphismes
de H sur $\PP^{n-1}\times\PP^{n-1}$. On a donc
$\text{deg}_{\ell}(p_{\F}^{*}(\mathcal{O}_{\PP^{n-1}}(1)_{|\ell})=1$
et $a_{\F}=a_{\G}=2$.\\
\indent Le lieu exceptionnel E s'identifie donc au recollement de deux
copies de
$\PP_{\PP^{n-1}}(\mathcal{E})$ le long de $\PP^{n-1}\times\PP^{n-1}$ et
l'automorphisme de
$\PP^{n-1}\times\PP^{n-1}$ induit par les deux isomorphismes de H sur
$\PP^{n-1}\times\PP^{n-1}$ est de la forme $(x,y)\mapsto
(\alpha(y),\,\beta(x))$ o\`u $\alpha$ et $\beta$ sont deux automorphismes
de $\PP^{n-1}$.\\
\indent Il reste \`a v\'erifier que tout
automorphisme de $\PP^{n-1}\times\PP^{n-1}$ de la forme
$(x,y)\mapsto (\alpha(x),\,\beta(y))$ o\`u $\alpha$ et $\beta$ sont deux
automorphismes de $\PP^{n-1}$ est la restriction \`a
$\PP^{n-1}\times\PP^{n-1}$ d'un automorphisme de
$\PP_{\PP^{n-1}}(\mathcal{E})$ stabilisant $\PP^{n-1}\times\PP^{n-1}$.\\
\indent Fixons un automorphisme $\beta$ de $\PP^{n-1}$ et notons
encore $\beta$ l'automorphisme lin\'eaire de
$\textup{W}\otimes\mathcal{O}_{\PP^{n-1}}$ correspondant.
L'automorphisme de $\PP_{\PP^{n-1}}(\mathcal{E})$ au dessus de
$\PP^{n-1}$ induit par $\text{Id}\oplus\beta$ stabilise
$\PP^{n-1}\times\PP^{n-1}$ et sa restriction \`a
$\PP^{n-1}\times\PP^{n-1}$ est $(x,y)\mapsto (x,\,\beta(y))$. Soit
$\alpha$ un automorphisme de $\PP^{n-1}$. Fixons des isomorphismes
$\alpha^{*}\mathcal{O}_{\PP^{n-1}}(2)\simeq\mathcal{O}_{\PP^{n-1}}(2)$
et
$\alpha^{*}\textup{W}\otimes\mathcal{O}_{\PP^{n-1}}\simeq\textup{W}\otimes
\mathcal{O}_{\PP^{n-1}}$.
L'isomorphisme
$\alpha^{*}\mathcal{E}\simeq\mathcal{E}$ d\'eduit des deux
pr\'ec\'edents induit un automorphisme de
$\PP_{\PP^{n-1}}(\mathcal{E})$ au dessus de $\alpha$ qui stabilise
$\PP^{n-1}\times\PP^{n-1}$ et sa restriction \`a
$\PP^{n-1}\times\PP^{n-1}$ est de la forme
$(x,y)\mapsto (\alpha(x),\,\beta'(y))$ o\`u $\beta'$ est un
automorphisme de $\PP^{n-1}$.\qed
\lemme{On a $\textup{H}^{1}(\F,\T_{\F}\otimes\mathcal{O}_{\F}(-\F))=
\textup{H}^{1}(\G,\T_{\G}\otimes\mathcal{O}_{\G}(-\G))=0$
et $\textup{H}^{1}(\E,\T_{\X}\otimes\mathcal{O}_{\E}(-i\E))=0$ pour tout 
$i\ge 1$ si $n\ge 3$.}
\textit{D\'emonstration}.$-$Les faisceaux
$\text{R}^{j}{p_{\F}}_{*}(\mathcal{O}_{\F}(k))$
sont nuls pour $j\in\{1,2\}$ et $k\in\ZZ$. Il existe donc
un isomorphisme $\displaystyle{\textup{H}^{i}(\F,\,\mathcal{O}_{\F}(k)\otimes
p_{\F}^{*}\mathcal{O}_{\PP^{n-1}}(l))\simeq
\textup{H}^{i}(\PP^{n-1},\,\textup{S}^{k}(\mathcal{O}_{\PP^{n-1}}(2)
\oplus\textup{W}\otimes\mathcal{O}_{\PP^{n-1}})\otimes
\mathcal{O}_{\PP^{n-1}}(l))}$
pour $i\in\{1,2\}$ et  $k,\,l\in\ZZ$. Le groupe
$\textup{H}^{i}(\F,\,\mathcal{O}_{\F}(k)\otimes
p_{\F}^{*}\mathcal{O}_{\PP^{n-1}}(l))$ est donc nul
pour $k,\,l\in\ZZ$ si $i=1$ et pour $l>-3$ si $i=2$. Les groupes
$\textup{H}^{1}(\F,\,p_{\F}^{*}(\textup{T}_{\PP^{n-1}}(l))\otimes
\mathcal{O}_{\F}(k))$
et $\textup{H}^{1}(\PP^{n-1},\text{T}_{\PP^{n-1}}(l)\otimes
\textup{S}^{k}(\mathcal{O}_{\PP^{n-1}}(2)\oplus\textup{W}\otimes
\mathcal{O}_{\PP^{n-1}}))$
sont \'egalement isomorphes pour $k,\,l\in\ZZ$ et donc nuls pour
$l>-3$.\\
\indent La longue suite exacte de cohomologie d\'eduite de la suite
exacte courte
$$(0\longrightarrow\mathcal{O}_{\F}\longrightarrow
p_{\F}^{*}(\mathcal{O}_{\PP^{n-1}}(-2)\oplus\textup{W}^{*}\otimes
\mathcal{O}_{\PP^{n-1}})\otimes\mathcal{O}_{\F}(1)
\longrightarrow \text{T}_{\F/\PP^{n-1}}\longrightarrow 0)\otimes
\mathcal{O}_{\F}(k)\otimes
p_{\F}^{*}\mathcal{O}_{\PP^{n-1}}(l)$$
et les annulations
pr\'ec\'edentes donnent l'annulation du groupe
$\textup{H}^{1}(\F,\text{T}_{\F/\PP^{n-1}}\otimes\mathcal{O}_{\F}(k)\otimes
p_{\F}^{*}\mathcal{O}_{\PP^{n-1}}(l))$ pour $l>-3$ et, de m\^eme, la
longue suite
exacte de cohomologie d\'eduite de la suite exacte courte
$$(0\longrightarrow\text{T}_{\F/\PP^{n-1}}
\longrightarrow\text{T}_{\F}\longrightarrow
p_{\F}^{*}\text{T}_{\PP^{n-1}}\longrightarrow
0)\otimes\mathcal{O}_{\F}(k)\otimes
p_{\F}^{*}\mathcal{O}_{\PP^{n-1}}(l)$$
donne finalement l'annulation du groupe
$\textup{H}^{1}(\F,\text{T}_{\F}\otimes\mathcal{O}_{\F}(k)\otimes
p_{\F}^{*}\mathcal{O}_{\PP^{n-1}}(l))$
pour $k+l>-3$. Le fibr\'e $\mathcal{O}_{\F}(-\F)$ est isomorphe au
fibr\'e
$\mathcal{O}_{\F}(2)\otimes p_{\F}^{*}\mathcal{O}_{\PP^{n-1}}(-1)$
et les groupes $\textup{H}^{1}(\F,\T_{\F}\otimes\mathcal{O}_{\F}(-\F))$ et
$\textup{H}^{1}(\G,\T_{\G}\otimes\mathcal{O}_{\G}(-\G))$ sont donc nuls.\\
\indent La longue suite exacte de cohomologie d\'eduite de la suite
exacte courte
$$(0\longrightarrow\text{T}_{\F}\longrightarrow{\text{T}_{\X}}_{|\F}
\longrightarrow\text{N}_{\F/\X}=\mathcal{O}_{\F}(-2)\otimes
p_{\F}^{*}\mathcal{O}_{\PP^{n-1}}(1)
\longrightarrow 0)\otimes\mathcal{O}_{\F}(k)\otimes
p_{\F}^{*}\mathcal{O}_{\PP^{n-1}}(l)$$ donne l'annulation de
$\textup{H}^{1}(\F,{\text{T}_{\X}}_{|\F}\otimes\mathcal{O}_{\F}(k)\otimes
p_{\F}^{*}\mathcal{O}_{\PP^{n-1}}(l))$ pour $l>-3$. L'id\'eal de E
dans $\X$ est $\I_{\F}\I_{\G}=\I_{\F}\cap\I_{\G}$ et la longue suite
exacte
de cohomologie d\'eduite de la suite exacte courte
$$(0\longrightarrow
\mathcal{O}_{\G}(-\F)=\mathcal{O}_{\G}(-1)\otimes
p_{\G}^{*}\mathcal{O}_{\PP^{n-1}}(2)
\longrightarrow\mathcal{O}_{\E}
\longrightarrow\mathcal{O}_{\F}\longrightarrow
0)\otimes{\textup{T}_{\X}}_{|\textup{E}}\otimes\mathcal{O}_{\E}(-k\E)$$
donne finalement l'annulation du groupe
$\textup{H}^{1}(\E,{\text{T}_{\X}}_{|\E}\otimes\mathcal{O}_{\E}(-k\E))$ pour
$k>-3$
et en particulier pour $k>0$.\qed\\
\newline
\marque Soit $\zeta$ une racine primitive cubique de l'unit\'e et soit $<\zeta>$
le groupe cyclique engendr\'e par $\zeta$. Notons
$\V'=\CC^{2n}\diagup<\zeta>$ o\`u $<\zeta>$ agit sur $\CC^{2n}$ par la
formule\\
\centerline{$\zeta.(z_{1},\,z_{2}\,,\ldots,\,z_{2n-1},\,z_{2n})=
(\zeta z_{1},\,\zeta^{2}z_{2},\,\ldots,\,\zeta z_{2n-1},\,\zeta^{2}z_{2n})$} et
$\X'$ l'\'eclatement normalis\'e de l'id\'eal maximal de $\textup{o}'$
dans $\V'$ o\`u $\textup{o}'$ est le point singulier de $\V'$.
Notons $\E'=\F'\cup\G'$ le diviseur exceptionnel et $\E_{i}$ (resp.
$\E'_{i}$)
le $i^{i\text{\`e}me}$ voisinage infinit\'esimal de $\E$ (resp. $\E'$)
dans $\X$ (resp. $\X'$) pour $i\ge 1$.
Notons $\F_{i}$ et $\G_{i}$ (resp. $\F'_{i}$ et $\G'_{i}$)
les $i^{i\text{\`e}mes}$ voisinages infinit\'esimaux de $\F$ et $\G$
(resp. $\F'$ et $\G'$) dans $\X$ (resp. $\X'$).
Fixons un isomorphisme $\psi_{1}\,:\,\E_{1}\simeq\E_{1}'$
(\textit{voir} proposition 3.5) qui applique
$\F$ (resp. $\G$) sur $\F'$ (resp. $\G'$).
\prop{Il existe un isomorphisme
$\psi_{2}\,:\,\E_{2}\simeq\E_{2}'$
compatible avec l'isomorphisme $\psi_{1}\,:\,\E_{1}\simeq\E_{1}'$.}
\textit{D\'emonstration}.$-$L'obstruction \`a l'existence d'un morphisme
$\E_{2}\longrightarrow\X'$
\'etendant le morphisme $\E\simeq\E'\hookrightarrow \X'$
est un \'el\'ement de $\textup{H}^{1}(\E,\T_{\X'}\otimes\mathcal{O}_{\E}(-\E))=
\textup{H}^{1}(\E',\T_{\X'}\otimes\mathcal{O}_{\E'}(-\E'))$ \re{Gr71}{Expos\'e 3
corollaire 5.2}
nul par le lemme 3.6. Fixons un tel morphisme. Il se factorise
\`a travers $\E'_{2}$ et d\'etermine un morphisme
$\psi_{2}\,:\,\E_{2}\longrightarrow\E'_{2}$ et,
par restriction \`a $\F_{2}$ (resp. $\G_{2}$), un morphisme
$\psi_{\F_{2}}\,:\,\F_{2}\longrightarrow\F'_{2}$ (resp.
$\psi_{\G_{2}}\,:\,\G_{2}\longrightarrow\G'_{2}$) \'etendant
la restriction $\psi_{\F}$ (resp. $\psi_{\G}$) de $\psi$ \`a $\F$
(resp. $\G$). Montrons que $\psi_{2}$ est en fait un isomorphisme.\\
\indent L'obstruction \`a l'existence d'un morphisme
$\F_{2}\longrightarrow\F$
(resp. $\F'_{2}\longrightarrow\F'$)
\'etendant l'identit\'e $\F\longrightarrow\F$ (resp.
$\F'\longrightarrow\F'$) est un \'el\'ement de
$\textup{H}^{1}(\F,\T_{\F}\otimes\mathcal{O}_{\F}(-\F))$ (resp.
$\textup{H}^{1}(\F',\T_{\F'}\otimes\mathcal{O}_{\F'}(-\F'))$)
nul par le lemme 3.6. Le sch\'ema $\F_{2}$ (resp. $\F'_{2}$) est donc un
$\F$-sch\'ema
(resp. $\F'$-sch\'ema) dont le faisceau structural est isomorphe au
quotient
$S^{\bullet}(\mathcal{O}_{\F}(-\F))/\mathcal{O}_{\F}(-2\F)
S^{\bullet}(\mathcal{O}_{\F}(-\F))$
(resp.
$S^{\bullet}(\mathcal{O}_{\F'}(-\F'))/\mathcal{O}_{\F'}(-2\F')
S^{\bullet}(\mathcal{O}_{\F'}(-\F'))$).
Les sch\'emas  ${\F}_{2}$ et ${\F'}_{2}$ sont donc isomorphes.\\
\indent Notons $i$ et $j$ les injections respectivement de
$\mathcal{O}_{\G'}$ et $\mathcal{O}_{\F'}(-\F')$ dans
$\mathcal{O}_{\F'}\oplus\mathcal{O}_{\F'}(-\F')$ et notons $p$ et $q$
les projections sur $\mathcal{O}_{\F'}$ et $\mathcal{O}_{\F'}(-\F')$
respectivement. Remarquons enfin
l'\'egalit\'e ${\psi_{\F_{2}}}_{*}\mathcal{O}_{\F_{2}}=
{\psi_{\F}}_{*}\mathcal{O}_{\F}\oplus{\psi_{\F}}_{*}\mathcal{O}_{\F}(-\F)=
\mathcal{O}_{\F'}\oplus\mathcal{O}_{\F'}(-\F')$. Le morphisme d'anneaux
${\psi_{\F_{2}}}^{*}\,:\,\mathcal{O}_{\F'}\oplus\mathcal{O}_{\F'}(-\F')
\longrightarrow\mathcal{O}_{\F'}\oplus\mathcal{O}_{\F'}(-\F')$ est
d\'etermin\'e par la d\'erivation
$\text{D}_{\F}=q\circ\psi_{\F_{2}}^{*}\circ
i\in\text{Der}_{\CC}(\mathcal{O}_{\F'},\,\mathcal{O}_{\F'}(-\F'))$
et l'application $\mathcal{O}_{\F'}$-lin\'eaire
$p\circ\psi_{\F_{2}}^{*}\circ
j\in\textup{Hom}_{\F'}(\mathcal{O}_{\F'}(-\F'),\,\mathcal{O}_{\F'}(-\F'))=\CC$.
Nous noterons $\lambda_{\F}$ l'\'el\'ement de $\CC$ correspondant \`a
$p\circ\psi_{\F_{2}}^{*}\circ i$. Notons enfin
$\text{D}_{\G}\in\text{Der}_{\CC}(\mathcal{O}_{\G'},\,\mathcal{O}_{\G'}(-\G'))$
et
$\lambda_{\G}\in\textup{Hom}_{\G'}(\mathcal{O}_{\G'}(-\G'),\,\mathcal{O}_{\G'}(-\G'))$
les \'el\'ements correspondants associ\'es \`a $\psi_{\G_{2}}$.\\
\indent Le morphisme $\psi_{\F_{2}}$ est un hom\'eomorphisme et si
$\lambda_{\F}\neq 0$ alors c'est un isomorphisme de
sch\'emas. Supposons donc $\lambda_{\F}\neq 0$ et v\'erifions que
$\psi_{\F_{2}}$ induit un isomorphisme de $\F_{2}\cap\G_{2}$ sur
$\F'_{2}\cap\G'_{2}$. Soit $\U\subset\F'$ un ouvert affine non vide
d'anneau $\A$
au dessus duquel le fibr\'e $\mathcal{O}_{\F'}(-\F')$ est trivial. Les
anneaux $\mathcal{O}_{\F_{2}}(\U)$ et $\mathcal{O}_{\F'_{2}}(\U)$
sont isomorphes \`a l'anneau $A\oplus\varepsilon A$ avec
$\varepsilon^{2}=0$. Le morphisme $\psi_{\F_{2}}^{*}$ est
$a+\varepsilon b\mapsto a+\varepsilon(\text{D}_{\F}(a)+\lambda_{\F}b)$. Soit
$h+\varepsilon h_{1}$ un g\'en\'erateur de l'id\'eal de $\F'_{2}\cap
\G'_{2}$ dans $\F'_{2}$ sur $\U$. L'\'el\'ement $h$
est une \'equation du diviseur $\F'\cap\G'_{2}$ dans $\F'$. La
restriction de $\psi_{2}$
\`a $\F_{2}\cap\G_{2}$ se factorise \`a travers $\F'_{2}\cap\G'_{2}$ et
l'\'el\'ement $\psi_{\F_{2}}^{*}(h+\varepsilon h_{1})$ est donc
dans l'id\'eal de $\F_{2}\cap\G_{2}$ dans $\F_{2}$. Ce dernier est
engendr\'e par un \'el\'ement de la forme $uh+\varepsilon h_{3}$ o\`u
$u$ est inversible dans $A$. Ecrivons
$\psi_{\F_{2}}^{*}(h+\varepsilon h_{1})=h+\varepsilon
h_{2}$ avec $h_{2}=\text{D}_{\F}(h)+\lambda_{\F}h_{1}$.
Il existe donc une relation
$h+\varepsilon h_{2}=(a+\varepsilon b)(uh+\varepsilon h_{3})
=auh+\varepsilon(ah_{3}+buh)$ et en particulier $h=auh$.
L'\'el\'ement $a$ est inversible dans $A$ et $a+\varepsilon b$ est donc
inversible dans $A\oplus\varepsilon A$. L'image par $\psi_{\F_{2}}^{*}$
de l'id\'eal de $A\oplus\varepsilon A$ engendr\'e par $h+\varepsilon
h_{1}$ est l'id\'eal - $\psi_{\F_{2}}^{*}$ est un isomorphisme -
engendr\'e par $\psi_{\F_{2}}^{*}(h+\varepsilon h_{1})$ et donc
l'id\'eal de $\F_{2}\cap\G_{2}$ dans $\F_{2}$.\\
\indent Si $\lambda_{\F}=0$ alors le morphisme $\psi_{\F_{2}}$ se
factorise
\`a travers $\F'\hookrightarrow\F'_{2}$.
L'image sch\'ematique de $\F_{2}\cap\G_{2}$ est en particulier un
sous-sch\'ema ferm\'e de $\F'$ et le morphisme
induit par $\psi_{2}$ de $\F_{2}\cap\G_{2}$ vers $\F'_{2}\cap\G'_{2}$
n'est pas un isomorphisme.\\
\indent En particulier, si $\lambda_{\F}=0$ alors $\lambda_{\G}=0$
et l'image sch\'ematique
de $\F_{2}\cap\G_{2}$ est donc le sous-sch\'ema ferm\'e $\F'\cap
\G'$, autrement dit, il existe une r\'etraction de l'inclusion
$\E\subset\E_{2}$, ce qui est exclu par le lemme 3.9.\\
\indent On a donc finalement $\lambda_{\F}\neq 0$, $\lambda_{\G}\neq 0$
et $\psi_{2}$ est
bien l'isomorphisme cherch\'e.\qed
\lemmeref{De01}{lemme 2.8}{Soit $(\R,\,\mathfrak{m})$ un anneau local
noeth\'erien et soit $\I\subset\mathfrak{m}^{2}$ un id\'eal de
$\R$. S'il existe une section de la projection canonique
$\R\twoheadrightarrow\R/\I$ alors $\I=0$.}
\indent Le lemme suivant et le lemme 3.6 nous assurent que
l'isomorphisme $\psi_{2}$
s'\'etend en un isomorphisme $\widehat{\X}\simeq\widehat{\X}'$ des
compl\'et\'es formels de X et $\X'$ le long de E et $\E'$
respectivement. Les compl\'et\'es
${\widehat{\mathcal{O}}}_{\V,\textup{o}}$ et
${\widehat{\mathcal{O}}}_{\V',\textup{o'}}$
sont donc isomorphes \re{Gr66}{Chap. III th\'eor\`eme 4.1.5}. Le th\'eor\`eme
d'approximation d'Artin \re{Ar68}{corollaire 1.6} termine 
la preuve du th\'eor\`eme.
\lemmeref{Mo82}{lemme 3.33}{Soient $\X$ et $\X'$ deux espaces
analytiques complexes r\'eduits
et irr\'eductibles de m\^eme dimension et $\E$ et $\E'$ deux diviseurs
de Cartier effectifs sur $\X$ et $\X'$ respectivement.
Notons $\E_{i}$ (resp. $\E'_{i}$) le $i^{i\text{\`e}me}$ voisinage
infinit\'esimal de $\E$ (resp. $\E'$) dans
$\X$ (resp. $\X'$) pour $i\ge 1$. On suppose que $\X$ est lisse et qu'il
existe un isomorphisme
$\E_{i_{0}}\simeq\E'_{i_{0}}$ pour un $i_{0}\ge 2$. Si
$\textup{H}^{1}(\E,\T_{\X}\otimes\mathcal{O}_{\E}(-i\E))=0$ pour $i\ge i_{0}$
alors l'isomorphisme
$\E_{i_{0}}\simeq\E'_{i_{0}}$ s'\'etend en un isomorphisme
$\widehat{\X}\simeq\widehat{\X}'$ des
compl\'et\'es formels de $\X$ et $\X'$ le long de $\E$ et $\E'$
respectivement.}

\noindent St\'ephane \textsc{Druel}, \textsc{Institut Fourier}, UMR 5582
du CNRS,
Universit\'e Joseph Fourier, BP 74, 38402 Saint Martin d'H\`eres,
France.\\
e-mail : \texttt{druel@mozart.ujf-grenoble.fr}

\end{document}